\documentclass[11pt]{article}
\usepackage{amssymb}
\usepackage{color}
\usepackage{soul}
\usepackage{ulem}
\newcommand{\bC}{\mathbf{C}}
\newcommand{\bN}{\mathbf{N}}
\newcommand{\bR}{\mathbf{R}}

\newcommand{\bK}{\mathbf{K}}
\newcommand{\ord}{\mathop{\mathrm{ord}}}

\newcommand{\Lo}{{\cal L}_0}

\newcommand{\qed}{\rule{0.5em}{0.5em}}

\newlength{\szer}

\newcommand{\Teissr}[2]{%
\settowidth{\szer}{$\displaystyle\frac{#1}{#2}$}%
\setlength{\szer}{0.5\szer}%
\left\{\hspace{\szer}%
\raisebox{1.5ex}{\makebox[0pt]{$#1$}}%
\raisebox{-1.5ex}{\makebox[0pt]{$#2$}}%
\hspace{-1\szer}\rule[0.4ex]{2\szer}{0.13ex}%
\hspace{-2\szer}\rule[0.7ex]{2\szer}{0.13ex}%
\right\}%
}

\newtheorem{defi}{Definition}[section]

\newtheorem{teorema}[defi]{Theorem}

\newtheorem{lema}[defi]{Lemma}
\newtheorem{coro}[defi]{Corollary}

\newenvironment{proof}[1][Proof]{\textbf{#1.} }{\
\rule{0.5em}{0.5em}}

\begin{document}
\title{\L ojasiewicz exponents   and Farey sequences
\footnotetext{
     \noindent  \begin{minipage}[t]{4in}{\small 2010 {\it Mathematics Subject Classification:\/} Primary 13F25;
      Secondary 14B05, 32S10.\\
       Key words and phrases: \L ojasiewicz exponent, logarithmic distance, Newton diagram, Farey sequences.\\
       The first-named and second-named authors were partially supported by the Spanish Project
       MTM 2012-36917-C03-01.}
       \end{minipage}}}

\author{A. B. de Felipe, E. R.\ Garc\'{\i}a Barroso, J. Gwo\'zdziewicz and A. P\l oski}

\maketitle

\begin{abstract}
\noindent Let $I$ be an ideal of the ring of formal power series $\bK[[x,y]]$ with coefficients in an algebraically closed 
field $\bK$ of arbitrary characteristic.
Let $\Phi$ denote the set of all parametrizations $\varphi=(\varphi_1,\varphi_2)\in \bK[[t]]^2$, where $\varphi \neq (0,0)$ and  
$\varphi (0,0)=(0,0)$. The purpose of this paper is to investigate the invariant
\[
\Lo(I)=\sup_{\varphi \in \Phi}\left( \inf_{f\in I} \frac{\ord f \circ \varphi}{\ord \varphi}\right)
\]

\noindent called the {\it \L ojasiewicz exponent} of $I$. Our main result states that for the ideals $I$ of finite 
codimension the \L ojasiewicz exponent $\Lo(I)$ is a Farey number i.e. an integer or a rational number of the 
form $N+\frac{b}{a}$, where $a,b,N$ are integers such that $0<b<a<N$.
\end{abstract}

\section{Introduction}
\label{intro}
\noindent  
Let $\bK$ be an algebraically closed field of arbitrary characteristic. Let $t$ be a variable. 
A {\it parametrization} is a pair $\varphi(t)=(\varphi_1(t),\varphi_2(t))\in \bK[[t]]^2\setminus \{(0,0)\}$ such that $\varphi_1(0)=\varphi_2(0)=0$. 
We put $\ord \varphi=\inf\{\ord \varphi_1,\ord \varphi_2\}$, 
where $\ord \varphi_k$ stands for the {\it order} of the power series $\varphi_k=\varphi_k(t)$. 
For any ideal $I\subset \bK[[x,y]]$ we consider  the {\it \L ojasiewicz exponent} $\Lo(I)$ 
(see \cite{L-T}, \cite{P}, \cite{D}, \cite{dF}, \cite{B-R}) defined by the formula

\[
\Lo(I)=\sup_{\varphi \in \Phi}\left( \inf_{f\in I} \frac{\ord f \circ \varphi}{\ord \varphi}\right),
\]
 
\noindent where $\Phi$ stands  for the set of all parametrizations $\varphi=(\varphi_1,\varphi_2)$.
 
\medskip 
 
\noindent Note that  $\Lo(I)< +\infty$ if and only if $I$ is of finite codimension. 

\medskip

\noindent In the framework of the complex analytic geometry the notion of \L ojasiewicz exponent was introduced and studied by Lejeune-Jalabert and Teissier \cite{L-T}. They considered much more general notion including the  \L ojasiewicz exponent of holomorphic ideals in several variables. D'Angelo  \cite{D} defined this invariant independently and gave its 
applications to complex function theory on domains in $\bC^n$. Recently Cassou-Nogu\`es and Veys \cite{CN-V} introduced an algorithm to study ideals in $\bK[[x,y]]$ which enables us to compute  $\Lo(I)$ using a finite sequence of Newton diagrams.

\medskip

\noindent  Let $g\in \bK[[x,y]]$ be an irreducible power series. We put 

\[ \Lo(I,g)=\inf_{f\in I}\left\{\frac{\ord f \circ \varphi}{\ord \varphi}\right\},\]

\noindent where $\varphi$ is a parametrization such that $g\circ \varphi=0$. The notion does not depend on the 
choice of $\varphi$. If $\Lo(I)=\Lo(I,g)$ then we say that the  \L ojasiewicz exponent {\it is attained on the branch} $g=0$.

\begin{teorema}(\cite[Theorem 6]{B-R})
\label{Th1}
Let $I\subset \bK[[x,y]]$ be a proper ideal and let $f_1,\ldots,f_m$ be generators of $I$.  Then
there is an irreducible factor $g$ of the power series $f_1 \cdots f_m$ such that $\Lo(I)$ is attained on the branch $g=0$.
\end{teorema}

\noindent This result was proved by Ch\c{a}dzy\'nski and Krasi\'nski \cite[Theorem 3]{Ch-K} and independently by McNeal and N\'emethi \cite[Theorem 1.1]{M-N} for holomorphic ideals. The case of ideals in $\bK[[x,y]]$, where $\bK$ is of arbitrary characteristic is due to Brzostowski and Rodak \cite[Theorem 6]{B-R}. In Section \ref{section2} of this note we give a very short proof of it. Let us write down the following corollary to Theorem \ref{Th1}.

\begin{coro}
If $I\subset \bK[[x,y]]$ is of finite codimension then $\Lo(I)$ is  a rational number. 
\end{coro}

\noindent Our main result is

\begin{teorema} 
\label{Th2}
Let $I$ be an ideal of $\bK[[x,y]]$ of finite codimension. Then $\Lo(I)$ is a Farey number i.e.  $\Lo(I)$ is an integer or a 
rational number of the form $N+\frac{b}{a}$, where $N,a,b$ are integers such that $0<b<a<N$.
\end{teorema}

\noindent Theorem \ref{Th2} gives a positive answer to Question 1 of \cite{B-R}. It implies that 
the fractional parts of the \L ojasiewicz exponents $\Lo(I)$ form the Farey sequences of order $\lfloor \Lo(I) \rfloor$ (see \cite{H-W}), where $\lfloor {z} \rfloor$ denotes the integer part of the real number $z$.

\medskip

\noindent The proof of Theorem \ref{Th2} is given in Section \ref{section3}.

\medskip

\noindent The holomorphic version ($I$ is a holomorphic ideal generated by two elements) was proved in \cite[Theorem 3.4]{P}. Its proof does not extend to the case of arbitrary characteristic. 

\section{Proof of Theorem \ref{Th1}}
\label{section2}

\noindent For any $f,g\in \bK[[x,y]]$ we consider the intersection number $i_0(f,g)=\dim_{\bK} \bK[[x,y]]/(f,g)$, where 
$(f,g)$ is the ideal generated by $f$ and $g$ in $\bK[[x,y]]$. Let 
\[d(f,g)=\frac{i_0(f,g)}{\ord f \,\ord g}\]

\noindent  for irreducible $f,g\in \bK[[x,y]]$. Then $d(f,g)$  is a {\it logarithmic distance} on the set of all irreducible power series, that is
\begin{enumerate}
\item[(D1)] $d(f,f)=+\infty$, 
\item[(D2)] $d(f,g)=d(g,f)$, and
\item[(D3)] $d(f,g)\geq \inf\{d(f,h),d(g,h)\}$ for $f,g,h$ irreducible power series.
\end{enumerate}

\noindent Only Property (D3) is non-trivial (see \cite[Corollary 2.9]{GB-P}). 

\medskip

\noindent If $g\in \bK[[x,y]]$ is irreducible then there exists a parametrization $\psi^o\in \bK[[t]]^2$ such that  $g\circ \psi^0=0$ and $\ord \psi^0=\ord g$. Moreover, for any power series $f\in \bK[[x,y]]$ we have $i_0(f,g)=\ord (f\circ \psi^o)$. 
If $\psi$ is a parametrization such that $g\circ \psi=0$ then 
there exists $\tau\in\bK[[t]]$ of positive order such that $\psi=\psi^o\circ \tau$. The equality 
$\frac{\ord (f\circ \psi)}{\ord \psi} = \frac{\ord (f\circ \psi^o)}{\ord \psi^o}=\frac{i_0(f,g)}{\ord g}$ 
shows that the definition of $\Lo(I,g)$ is correct and  can be rewritten  as follows
$\Lo(I,g)=\displaystyle \inf_{f \in I}\frac{i_0(f,g)}{\ord g}$.

\noindent If $\varphi$ is a parametrization then there exists an irreducible power series $g\in \bK[[x,y]]$
such that $g\circ \varphi=0$. This shows that 
\[
\Lo(I)=\sup \left \{\Lo(I,g)\;:\; \hbox{\rm $g$ is irreducible} \right\}.
\]

\medskip

\noindent  If $I=(f_1,\ldots,f_m)$  then
\begin{equation}
\label{eq:1}
\Lo(I,g)=\displaystyle \inf_{1\leq k\leq m} \frac{i_0(f_k,g)}{\ord g}.
\end{equation}

\medskip

\begin{lema}
\label{L1}
Let $I=(f_1,\ldots,f_m)$
and let $\prod_i f_{i}=\prod_j h_j$ with $h_j\in\bK[[x,y]]$ irreducible.
Let $g\in \bK[[x,y]]$ be an irreducible power series. 
Take an index~$k$ such that 
$d(g,h_{k})=\sup_{j}\{d(g,h_j)\}$. Then $\Lo(I,g)\leq \Lo(I,h_{k}).$
\end{lema}

\noindent \begin{proof}
\noindent Let us denote $h_k$ by $h$.
Then $d(g,h)\geq d(g,h_j)$ for any index $j$. After (D3) we get 
$d(h,h_j)\geq \inf \{d(g,h),d(g,h_j)\}=d(g,h_j)$. 
Therefore for any $j$ we have 
$\frac{i_0(g,h_j)}{\ord g}\leq \frac{i_0(h,h_j)}{\ord h}$ 
and consequently for any $i\in \{1,\ldots,m\}$ we get
$\frac{i_0(g,f_{i})}{\ord g}\leq \frac{i_0(h,f_i)}{\ord h},$ 
which implies $\Lo(I,g)\leq \Lo(I,h)$. 
\end{proof}

\medskip

\noindent Now, we can prove Theorem \ref{Th1}. 

\medskip

\noindent {\bf Proof of Theorem \ref{Th1}.} 
We keep the notations of Lemma \ref{L1}.
Fix an irreducible power series $g$. 
Then $\Lo(I,g)\leq \Lo(I,h_k)$. 
Hence  $\Lo(I)\leq \sup_j \{\Lo(I,h_j)\}$. The inverse inequality is obvious. 
Therefore $\Lo(I)=\sup_j \{\Lo(I,h_j)\}$, which proves Theorem 1.1. \qed

\section{Proof of Theorem \ref{Th2}}
\label{section3}

\noindent Let $f=\sum c_{\alpha \beta} x^{\alpha}y^{\beta}\in \bK[[x,y]]$. The {\it Newton diagram} $\Delta(f)$ of $f$ is by definition the convex hull of the set $\{(\alpha,\beta)\in \bN^2\;:\; c_{\alpha \beta}\neq 0\}+ \bR^2_{\geq 0}$. We use Teissier's notation (\cite[p. 846]{L-T}) denoting by $\Teissr{b}{a}$ the Newton diagram of $y^a+x^b$, for $a,b>0$. The following properties of Newton diagrams are well-known
\begin{enumerate}
\item[(N1)] for generic $c_1,\ldots,c_m$: $\Delta(\sum_{i=1}^m c_if_i)$ is the convex hull of the set  $\bigcup_{i=1}^m \Delta(f_i)$, 
\item[(N2)] if  $f=0$ is a branch different from the axes then $\Delta(f)=\Teissr{i_0(f,y)}{i_0(f,x)}$,
\item[(N3)] if $\Delta(f_1)=\Teissr{b_1}{a_1}$ and $\Delta(f_2)=\Teissr{b_2}{a_2}$ then $i_0(f_1,f_2)\geq \min \{a_1b_2, a_2b_1\} $, with equality if $a_1b_2\neq a_2b_1$.
\end{enumerate}

\medskip

\noindent  Property (N1) is a consequence of the definition of $\Delta(f)$. For Property (N2) see \cite[Proposition 4.2]{Ploski2013}. Property (N3) follows from \cite[Proposition 3.13, Proposition 3.8 (v)]{Ploski2013}.

\medskip

\noindent 
Let $I$ be an ideal of $\bK[[x,y]]$ with a finite \L{}ojasiewicz exponent.
Put $l=\Lo(I)$. Consider the set of ideals $J \subset \bK[[x,y]]$  such that $\Lo(J)=l$ 
and let $M$ be a maximal  element of this set (with respect to the inclusion).
 Set $\ord M=\inf \{\ord f\;:\; f\in M\}$. 
 Observe that replacing any finite sequence of generators of $M$ 
 by their general linear combinations we obtain  generators of the same order, 
 equal to $\ord M$.

\begin{lema}
\label{L2}
If $f_1,\ldots, f_m$ is a sequence of generators of $M$ of the same order then there exists  
$ k\in \{1, ..., m\}$, such that $f_k$ is irreducible and  $\Lo(M,f_k)=\Lo(M)$.
\end{lema}
\noindent \begin{proof}
\noindent Let $f_1,\ldots,f_m$ be a sequence of generators of $M$ of the same order. By Theorem \ref{Th1} the \L ojasiewicz exponent of $M$ is attained on an irreducible factor $h$ of the product $f_1\cdots f_m$.

\medskip
\noindent Let $\overline M$ be the ideal generated by $f_1,\ldots,f_m$ and $h$. Since $M\subset \overline M$ we get $\Lo(M) \geq \Lo(\overline M)$. On the other hand $\Lo(\overline M)\geq \Lo(\overline M,h)= \Lo(M,h)=\Lo(M)$,  which implies
$\Lo(\overline M)=\Lo(M)$. By the maximality of $M$ we get $h\in M$. Let $k \in \{1,\ldots,m\}$ be an index such that $h$ divides $f_k$. Then $\ord h\leq \ord f_k=\ord M\leq \ord h$, hence $\ord h=\ord f_h$ and $f_k=h\cdot u$, where $u(0,0)\neq 0$ which implies that $f_k$ is irreducible.
\end{proof}

\medskip
\noindent Let us pass to the proof of Theorem \ref{Th2}. 

\medskip

\noindent {\bf Proof of Theorem \ref{Th2}.}
 We keep the notation and assumptions introduced above. 
It suffices to prove that $\Lo(M)$ is an integer or  $\Lo(M)=N+\frac{b}{a}$, where $0<b<a<N$.

\medskip
\noindent 
By Lemma \ref{L2} there exists an irreducible power series $h\in M$ 
of order $\ord h=\ord M$ such that $\Lo(M,h)=\Lo(M)$. 

\medskip
\noindent 
If $\ord h=1$ then $\Lo(M)$ is an integer.

\medskip
\noindent 
If $\ord h >1$ then changing the system of coordinates if necessary,  we may assume that 
$\ord h(0,y)< \ord h(x,0)$.
Let $a=\ord h(0,y)$ and $c=\ord h(x,0)$. 
Then $\Delta(h)=\Teissr{c}{a}$, where $a=\ord h=\ord M$.

\medskip
\noindent  
Replacing any finite sequence of generators of $M$ by a sequence of their linear generic combinations we get a sequence $f_1,\ldots,f_m$ of generators of the same order such that $\Delta(f_1)=\cdots=\Delta(f_m)$. 
Let $\Delta$ be their common Newton diagram. 

\medskip
\noindent 
Since $h\in M$, we have $h=a_1f_1+\cdots+a_mf_m$, where $a_i\in \bK[[x,y]]$.
Substituting $x=0$ we get 
$\ord M =\ord h(0,y)\geq \min_i\{\ord f_{i}(0,y)\} \geq \ord M$.
 Hence the diagram $\Delta$ intersects the vertical axis  at $(0,a)$. 
By Lemma \ref{L2} at least one of $f_1, \ldots,f_m$ is irreducible. 
This implies that $\Delta$ has only one compact face. 
Since $\ord f_i = a$, we have $\Delta=\Teissr{d}{a}$, where $d \geq a$. 

\medskip
\noindent 
By (\ref{eq:1}) there is $k\in\{1,\ldots,m\}$ such $\Lo(M)=\frac{i_0(f_k,h)}{\ord h}$.

\medskip
\noindent 
If $d=a$ then by  (N3) we get $i_0(f_k,h) = \min\{ac,a^2\}=a^2$. In this case 
$\Lo(M)=a$. 

\medskip
\noindent 
If $d>a$ then by  (N3) we get $i_0(f_k,h)\geq   \min \{ac,ad\}=a\min\{c,d\}\geq a(a+1)$. 
Write $i_0(f_k,h)=aN+b$, where $0\leq b<a$. 
Dividing this equality by $a$ and taking integer parts we get
$N=\left \lfloor \frac{i_0(f_{i_0},h)}{a}\right \rfloor \geq \frac{a(a+1)}{a}=a+1$. 
Therefore $\Lo(M)=N+\frac{b}{a}$, where $0\leq b<a<N$, which completes the proof.
\qed

\medskip
\noindent

\noindent {\small Ana Bel\'en de Felipe\\
Institut de Math\'ematiques de Jussieu-Paris Rive Gauche,\\
 B\^atiment Sophie Germain, case 7012,\\
75205 Paris Cedex 13, France.\\
ana.de-felipe@imj-prg.fr}

\medskip

\noindent {\small  Evelia Rosa Garc\'{\i}a Barroso\\
Departamento de Matem\'aticas, Estad\'{\i}stica e I.O. \\
Secci\'on de Matem\'aticas, Universidad de La Laguna\\
Apartado de Correos 456\\
38200 La Laguna, Tenerife, Espa\~na\\
e-mail: ergarcia@ull.es}

\medskip

\noindent {\small   Janusz Gwo\'zdziewicz\\
Institute of Mathematics\\
Pedagogical University of Cracow\\
Podchor\c a{\accent95 z}ych 2\\
PL-30-084 Cracow, Poland\\
e-mail: gwozdziewicz@up.krakow.pl}

\medskip

\noindent {\small Arkadiusz P\l oski\\
Department of Mathematics and Physics\\
Kielce University of Technology\\
Al. 1000 L PP7\\
25-314 Kielce, Poland\\
e-mail: matap@tu.kielce.pl}

\end{document}